\documentclass[12pt]{article}
\usepackage{latexsym, amssymb}
    \title{{\bf  Conformal-field-theoretic analogues of 
codes and lattices}}
    \author{Yi-Zhi Huang}
    \date{}
    \begin{document}
    \bibliographystyle{alpha}
    \maketitle

\newtheorem{thm}{Theorem}[section]
\newtheorem{defn}[thm]{Definition}
\newtheorem{prop}[thm]{Proposition}
\newtheorem{cor}[thm]{Corollary}
\newtheorem{lemma}[thm]{Lemma}
\newtheorem{rema}[thm]{Remark}
\newtheorem{app}[thm]{Application}
\newtheorem{prob}[thm]{Problem}
\newtheorem{conv}[thm]{Convention}
\newtheorem{conj}[thm]{Conjecture}
\newtheorem{expl}[thm]{Example}

\newcommand{\halmos}{\rule{1ex}{1.4ex}}
\newcommand{\pfbox}{\hspace*{\fill}\mbox{$\halmos$}}

\newcommand{\nn}{\nonumber \\}

 \newcommand{\res}{\mbox{\rm Res}}
\renewcommand{\hom}{\mbox{\rm Hom}}
 \newcommand{\pf}{{\it Proof.}\hspace{2ex}}
 \newcommand{\epf}{\hspace*{\fill}\mbox{$\halmos$}}
 \newcommand{\epfv}{\hspace*{\fill}\mbox{$\halmos$}\vspace{1em}}
 \newcommand{\epfe}{\hspace{2em}\halmos}
\newcommand{\nord}{\mbox{\scriptsize ${\circ\atop\circ}$}}
\newcommand{\wt}{\mbox{\rm wt}\ }
\newcommand{\swt}{\mbox{\rm {\scriptsize wt}}\ }
\newcommand{\clr}{\mbox{\rm clr}\ }
\newcommand{\lbar}{\bigg\vert}

\begin{abstract} 
We introduce and study completely-extendable conformal
intertwining algebras.  Based on 
results obtained in other papers, various examples are given. 
Duals of these algebras are
constructed and nondegenerate such 
algebras are defined. We prove that the double dual of such a
nondegenrate  algebra is equal to itself. 
 We explain using a table that these
nondegenerate  algebras are
the correct conformal-field-theoretic analogues of linear binary codes
and nondegenerate rational lattices.  
\end{abstract}

\renewcommand{\theequation}{\thesection.\arabic{equation}}
\renewcommand{\thethm}{\thesection.\arabic{thm}}
\setcounter{equation}{0}
\setcounter{thm}{0}
\setcounter{section}{-1}

\section{Introduction}

Codes and lattices play an important role in the study of many
mathematical problems and, in particular, they are essential in the study
of finite sporadic simple groups (see \cite{CN} for details). 
In \cite{B}, Borcherds introduced a notion of vertex algebra and
announced that the moonshine module constructed by 
Frenkel, Lepowsky and Meurman in \cite{FLM1} 
has a natural structure of vertex algebra. 
In \cite{FLM}, Frenkel, Lepowsky and Meurman defined
vertex operator algebras, which are vertex algebras equipped
with additional data and satisfying
additional axioms, and proved that the moonshine module
indeed has a structure of a vertex operator algebra (see also 
\cite{DGM} and \cite{H2} for different proofs).
In this major work, Frenkel, Lepowsky and Meurman
discovered that vertex operator algebras  have many properties 
analogous to properties of doubly-even codes and even
lattices (see \cite{FLM} for details). They emphasized that 
linear binary codes, nondegenerate 
rational lattices and conformal field theories should be viewed 
as three stages in  a hierarchy and the results in one stage should
have corresponding results in other stages.

For the Golay code and Leech
lattice, there are uniqueness results which state that any code
or lattice satisfying three suitable conditions must be
the Golay code or Leech
lattice (see \cite{C} and \cite{CN}). Based on the analogy above, especially
the analogy among the Golay code, Leech lattice and moonshine module
vertex operator algebra, Frenkel, Lepowsky and Meurman 
conjectured in \cite{FLM}
that the moonshine module should also have an analogous uniqueness
property.
In \cite{G}, Goddard discussed this analogy in more detail and
pointed out  that one needs to find suitable
nonmeromorphic 
generalizations of vertex
(operator) algebras in order to establish a precise correspondence.  In
particular, one needs a functor on the category
of these generalizations 
corresponding to the dual functors for codes and lattices.

Several nonmeromorphic generalizations of the notion of vertex operator
algebra were first given by 
by Frenkel, Feingold and Ries \cite{FFR} and
Dong and Lepowsky \cite{DL1} \cite{DL2}. The most general one among
these notions is the notion of abelian intertwining algebra in
\cite{DL1} and \cite{DL2}. The Verlinde algebras associated to these
generalizations are abelian group algebras and involve only
one-dimensional representations of the braid groups. A 
``non-abelian'' notion of intertwining operator algebra was introduced in
\cite{H4} and studied in \cite{H5} and \cite{H6}. Various
examples of intertwining operator algebras were constructed by
Lepowsky, Milas and the author
in \cite{H2}, \cite{H3}, \cite{HL4}, \cite{HM1} and \cite{HM2} 
using the results obtained in \cite{HL1}, \cite{HL2}, \cite{HL3}
and \cite{H1}. 

In the present paper, we introduce a notion of completely-extendable
conformal intertwining algebra.  The examples of intertwining operator
algebras constructed in \cite{H2}, \cite{H3}, \cite{HL4}, \cite{HM1}
and \cite{HM2} are examples of completely-extendable conformal
intertwining algebras (actually they are ``complete'').  These
examples are constructed using a main theorem in \cite{H1} proving the
associativity of intertwining operators. This theorem has a
generalization to conformal intertwining algebras. We state this
theorem in the present paper but its proof will be given in another
paper \cite{H8}.  We also give some other examples, including the
examples constructed from free bosons.  Duals of completely-extendable
conformal intertwining algebras are constructed and the notion of
nondegenerate completely-extendable intertwining operator algebras is
introduced in terms of duals. We prove that the double dual of such an
algebra is equal to itself.  A table filling some blanks in a table
given by Goddard in \cite{G} and containing more items than the one in
\cite{G} is given to demonstrate that nondegenerate
completely-extendable intertwining operator algebras are the correct
conformal-field-theoretic analogues of codes and lattices.  We hope
that the analogues and correspondence given in this paper will lead to
at least a strategy to the solution of the uniqueness conjecture of
Frenkel, Lepowsky and Meurman for the moonshine module vertex operator
algebra. A more conceptual formulation of this correspondence in terms
of the language of categories and related results on the correspondence
will be given in a future paper.

The present paper is organized as follows: We introduce the notion 
of locally-grading-restricted 
conformal intertwining (super)algebra and state the basic properties of 
these algebras in Section 2. The notion of completely-extendable
conformal intertwining (super)algebra is introduced in Section 3. 
Examples of these algebras
are also given in this section. The dual of 
a completely-extendable conformal
intertwining  algebra is constructed in Section 4. 
In the same section, the notion of nondegenerate 
completely-extendable
conformal intertwining algebra is 
introduced and it is proved that the double
dual of such an algebra is equal to itself. 
In Section 5, we discuss the analogy between codes, lattices and 
nondegenerate 
completely-extendable conformal
intertwining algebras.

\paragraph{Acknowledgment} A large part of the research in the present paper
was finished when I visited University of Virginia in the 
year of 2000. I am grateful to the Department of Mathematics 
in University of Virginia and  Brian Parshall for the hospitality. 
This research is supported in part 
by NSF grant DMS-0070800.

\renewcommand{\theequation}{\thesection.\arabic{equation}}
\renewcommand{\thethm}{\thesection.\arabic{thm}}
\setcounter{equation}{0} \setcounter{thm}{0}

\section{Conformal intertwining algebras}

In this section, we relax the grading-restriction conditions 
for intertwining operator 
(super)algebras to introduce  locally-grading-restricted 
conformal intertwining
(super)algebras. We also give two basic properties of these algebras.  Since
all these properties are trivial generalizations of the properties of
intertwining operator algebras, we omit the proofs here. For
details on intertwining operator algebras, see \cite{H4}, \cite{H5}
and \cite{H6}.  We assume that the reader is familiar with the basic
notions of vertex operator algebra, module and intertwining operator,
as presented in \cite{FLM} and \cite{FHL}.

First we need the following notion of 
locally-grading-restricted conformal vertex (super)algebra:

\begin{defn}
{\rm A {\it conformal vertex algebra} of central charge $c$
is a $\mathbb{Z}$-graded vector space
equipped with a {\it vertex operator map} $Y: V\otimes V\to V[[x, x^{-1}]]$
and two distinguished elements $\mathbf{1}$ (the {\it vacuum})
and $\omega$ (the {\it Virasoro element} or {\it conformal element})
satisfying all the axioms for vertex operator algebras 
of central charge $c$
except for the two grading-restriction axioms $\dim V_{(n)}<\infty$
for $n\in \mathbb{Z}$ and  $V_{(n)}=0$ when the real part of 
$n$ is sufficiently small. 
A conformal vertex 
algebra is said to be {\it locally grading-restricted} if (i)
for any $m>0$, $v_{1}, \dots, v_{m}\in V$, there exists
$r\in \mathbb{Z}$ such that
the coefficients of the series 
$Y(v_{1}, x_{1})\dots Y(v_{m-1}, x_{m-1})v_{m}$
is in $\coprod_{n>r}V_{(n)}$ and (ii)
for any element of the conformal vertex
algebra, the module $W=\coprod_{n\in \mathbb{Z}}W_{(n)}$ for the
Virasoro algebra generated by this element satisfies the
grading-restriction conditions, that is, $\dim W_{(n)}<\infty$ for $n\in
\mathbb{Z}$ and $W_{(n)}=0$ when the real part of $n$ is
sufficiently small. {\it Modules}
 for a locally-grading-restricted conformal
vertex algebra are defined in the same way as modules
for a vertex operator algebra except that
they are  required to be only locally grading-restricted in the 
sense above. {\it Intertwining operators} 
for a locally-grading-restricted conformal vertex algebra are
defined in the obvious way. All the concepts for vertex operator algebras,
for example, simple vertex operator algebras and fusion 
rules, can be generalized to locally-grading-restricted 
conformal vertex algebras without any difficulty.

{\it Conformal vertex superalgebras} are defined similarly except that the 
underlying vector space has an additional $\mathbb{Z}_{2}$ grading
in the commutativity, the Jacobi identity or the skew-symmetry:
When both elements are odd,
there is an extra minus sign in the term in which the 
order of the two elements is changed. All the other notions for
conformal vertex superalgebras are defined in the obvious way.}
\end{defn}

\begin{rema}
{\rm Condition (i) in the definition above,  was first stated in 
\cite{DL2} for abelian intertwining algebras.
It guarantees that the convergence and rationality 
of the matrix elements
of products and 
iterates of vertex operators. Condition (ii) makes sure that all
the results involving the Virasoro operators (for 
example, the geometry of vertex operator algebras in \cite{H0} 
and \cite{H4.5})
still hold for these
algebras.}
\end{rema}

We now define locally-grading-restricted conformal
intertwining (super)algebra.

\begin{defn}\label{ioa}
{\rm A {\it locally-grading-restricted conformal
intertwining algebra of central charge 
$c\in \mathbb{C}$} consists of the following data:

\begin{enumerate}

\item A vector space
$$W=\coprod_{a\in
\mathcal{A}}W^{a}$$
graded
by a set $\mathcal{A}$
containing a special element $e$.

\item A structure of locally-grading-restricted conformal vertex algebra
 of central charge $c$ 
on $W^{e}$, and a $W^{e}$-module 
structure on $W^{a}$ for 
each $a\in \mathcal{A}$.

\item A subspace $\mathcal{V}_{a_{1}a_{2}}^{a_{3}}$ of 
the space of all intertwining operators of type 
${W^{a_{3}}\choose W^{a_{1}}W^{a_{2}}}$ for  each triple
$a_{1}, a_{2}, a_{3}\in \mathcal{A}$.

\end{enumerate}

\noindent These data satisfy the
following axioms:

\begin{enumerate}

\item The $W^{e}$-module structure on $W^{e}$ is the adjoint module
structure. For any $a\in \mathcal{A}$, the space $\mathcal{V}_{ea}^{a}$ is the
one-dimensional vector space spanned by the vertex operators defining the
$W^{e}$-module structure on $W^{a}$.
For any $a_{1}, a_{2}\in \mathcal{A}$ such that
$a_{1}\ne a_{2}$, $\mathcal{V}_{ea_{1}}^{a_{2}}=0$.

\item {\it Convergence properties}: For any $m\in \mathbb{Z}_{+}$,
$a_{i}, b_{0}, b_{i}, \in \mathcal{A}$, $w_{(a_{i})}
\in W^{a_{i}}$, $\mathcal{Y}_{i}\in \mathcal{
V}_{a_{i}\;b_{i}}^{b_{i-1}}$, $i=1, \dots, m$, $w_{(b_{0})}'
\in (W^{b_{0}})'$ and 
$w_{(b_{m})}\in W^{b_{m}}$, the series
\begin{equation}\label{conv-pr}
\langle w_{(b_{0})}', \mathcal{Y}_{1}(w_{(a_{1})}, x_{1})
\cdots\mathcal{Y}_{m}(w_{(a_{m})},
x_{m})w_{(b_{m})}\rangle_{W^{\mu_{1}}}|_{x^{n}_{i}=e^{n\log z_{i}},
i=1, \dots, m, n\in \mathbb{C}}
\end{equation}
is absolutely convergent when $|z_{1}|>\cdots >|z_{m}|>0$, and 
consequently (using also other axioms) for 
any $\mathcal{Y}_{1}\in \mathcal{
V}_{a_{1}a_{2}}^{a_{5}}$ and $\mathcal{Y}_{2}\in \mathcal{
V}_{a_{5}a_{3}}^{a_{4}}$, the series
\begin{equation}\label{conv-it}
\langle w_{(a_{4})}', \mathcal{Y}_{2}(\mathcal{Y}_{1}(w_{(a_{1})}, 
x_{0})w_{(a_{2})},
x_{2})w_{(a_{3})}\rangle_{W^{a_{4}}}
|_{x^{n}_{0}=e^{n\log (z_{1}-z_{2})},
x^{n}_{2}=e^{n\log z_{2}}, n\in \mathbb{C}}
\end{equation}
is absolutely convergent when
$|z_{2}|>|z_{1}-z_{2}|>0$.

\item {\it Associativity}:  For any $a_{1}, a_{2}, 
a_{3}, a_{4}, a_{5}\in \mathcal{A}$, 
any $\mathcal{Y}_{1}\in \mathcal{
V}_{a_{1}a_{5}}^{a_{4}}$ and $\mathcal{Y}_{2}\in\mathcal{
V}_{a_{2}a_{3}}^{a_{5}}$, there exist 
$\mathcal{Y}^{a, i}_{3}
\in \mathcal{
V}_{a_{1}a_{2}}^{a}$ and $\mathcal{Y}^{a,i}_{4}\in \mathcal{
V}_{aa_{3}}^{a_{4}}$ for $a\in \mathcal{A}$ and $i=1, \dots, k$
such that only finitely many of $\mathcal{Y}^{a, i}_{3}$ and 
$\mathcal{Y}^{a,i}_{4}$, 
$a\in \mathcal{A}$, $i=1, \dots, k$, are nonzero, and 
for any $w_{(a_{i})}\in W^{a_{i}}$, $i=1, 2, 3$, and $w_{(a_{4})}'
\in W'_{a_{4}}$,
the (multivalued) analytic function 
\begin{equation}\label{prod}
\langle w_{(a_{4})}', 
\mathcal{Y}_{1}(w_{(a_{1})}, x_{1})\mathcal{Y}_{2}(w_{(a_{2})}, 
x_{2})w_{(a_{3})}\rangle_{W}|_{x_{1}=z_{1},
x_{2}=z_{2}}
\end{equation}
defined in the region
$|z_{1}|>|z_{2}|>0$ 
and the (multivalued) analytic function
\begin{equation}\label{iter}
\sum_{a\in \mathcal{A}}\sum_{i=1}^{k}
\langle w_{(a_{4})}', \mathcal{Y}^{a, i}_{4}
(\mathcal{Y}^{a, i}_{3}(w_{(a_{1})},
x_{0})w_{(a_{2})}, x_{2})w_{(a_{3})}\rangle_{W^{a_{4}}}
|_{x_{0}=z_{1}-z_{2},
x_{2}=z_{2}}
\end{equation} 
defined in the region
$|z_{2}|>|z_{1}-z_{2}|>0$ are equal in the intersection
$|z_{1}|> |z_{2}|>|z_{1}-z_{2}|>0$.

\item {\it Skew-symmetry}:  For any $a_{1}, a_{2}, 
a_{3}\in \mathcal{A}$, the restriction of  $\Omega_{-1}$ to
$\mathcal{V}_{a_{1}a_{2}}^{a_{3}}$ is an isomorphism from
$\mathcal{V}_{a_{1}a_{2}}^{a_{3}}$ to 
$\mathcal{V}_{a_{2}a_{1}}^{a_{3}}$. (Here
$\Omega_{-1}$ is the linear map
from the space of all intertwining operators of type
${a_{3}\choose a_{1}a_{2}}$ to the space of 
all intertwining operators of type
${a_{3}\choose a_{2}a_{1}}$ defined by
formula
$$
\Omega_{-1}(\mathcal{Y})(w_{(a_{2})},x)w_{(a_{1})} = e^{xL(-1)}
\mathcal{Y}(w_{(a_{1})},e^{-\pi i}x)w_{(a_{2})}.
$$
for  $w_{(a_{1})}\in W^{a_{1}}$ and $w_{(a_{2})}\in W^{a_{2}}$.
Cf. \cite{HL2}.)

\end{enumerate}

A locally-grading-restricted conformal 
intertwining algebra is said to be {\it meromorphic} if 
intertwining operators are all Laurent series and is said to be 
{\it $\mathbb{Z}$-graded} if weights of elements of $W$ are 
all integers.

An {\it intertwining operator algebra} is a 
locally-grading-restricted conformal 
intertwining algebra such that
$W^{e}$ is a vertex operator algebra,
$W^{a}$, $a\in \mathcal{A}$, are (grading-restricted) modules for the
vertex operator algebra $W^{e}$ and $\mathcal{A}$ is finite.
Similarly we have {\it meromorphic} intertwining operator algebras
and {\it $\mathbb{Z}$-graded} intertwining operator algebras.

{\it Locally-grading-restricted conformal intertwining superalgebras}
and {\it intertwining operator superalgebras} are defined similarly
in the obvious way.}
\end{defn}

The locally-grading-restricted conformal 
intertwining (super)algebra or intertwining operator (super)algebra
just defined above is denoted by 
$$(W, 
\mathcal{A}, \{\mathcal{V}_{a_{1}a_{2}}^{a_{3}}\}, {\bf 1}, 
\omega)$$ 
or simply $W$. For simplicity, below we shall often use the abbreviation 
``l.g.r.'' for the word ``locally-grading-restricted.''
Clearly,  $\mathbb{Z}$-graded l.g.r. conformal 
intertwining algebras are meromorphic.

\begin{rema}\label{super1}
{\rm Actually, l.g.r. conformal intertwining superalgebras
can be viewed as conformal intertwining algebras. We shall not 
discuss this in detail in this paper. See \cite{H8}. See also 
Remark \ref{super2} below.}
\end{rema}

Now we give two basic properties of l.g.r. conformal 
intertwining algebras. The proofs of these properties are completely the 
same as the proofs of the same properties for intertwining operator
algebras and will be omitted. The similar properties 
with possible extra signs 
for l.g.r. conformal intertwining superalgebras also hold. But we 
omit the discussions of them in this paper.

The following commutativity generalizes the commutativity 
for intertwining operator algebras proved in \cite{H3}:

\begin{prop}\label{commu}
Let $(W, 
\mathcal{A}, \{\mathcal{V}_{a_{1}a_{2}}^{a_{3}}\}, {\bf 1}, \omega)$ 
be a l.g.r. conformal 
intertwining algebra. Then we have the following
{\it commutativity}: For any $a_{1}, \dots, a_{5}\in \mathcal{A}$,
any $\mathcal{Y}_{1}\in \mathcal{
V}_{a_{1}a_{5}}^{a_{4}}$ and $\mathcal{Y}_{2}\in\mathcal{
V}_{a_{2}a_{3}}^{a_{5}}$, there exist 
$\mathcal{Y}^{a,i}_{5}
\in \mathcal{
V}_{a_{2}a}^{a_{4}}$ and $\mathcal{Y}^{a, i}_{6}\in \mathcal{
V}_{a_{2}a_{3}}^{a}$ for $a\in \mathcal{A}$ and $i=1, \dots, k$
such that that only finitely many of $\mathcal{Y}^{a, i}_{5}$ and 
$\mathcal{Y}^{a,i}_{6}$, 
$a\in \mathcal{A}$, $i=1, \dots, k$, are nonzero, and 
for any $w_{(a_{i})}\in W^{a_{i}}$, $i=1, 2, 3$, and $w_{(a_{4})}'
\in W'_{a_{4}}$,
the (multivalued) analytic function 
(\ref{prod})
defined in the region
$|z_{1}|>|z_{2}|>0$ 
and the (multivalued) analytic function
\begin{equation}\label{op-prod}
\sum_{a\in \mathcal{A}}\sum_{i=1}^{k}\langle w_{(a_{4})}', 
\mathcal{Y}^{a, i}_{5}(w_{(a_{2})}, x_{2})\mathcal{Y}^{a, i}_{6}(w_{(a_{1})}, 
x_{1})w_{(a_{3})}\rangle_{W^{a_{4}}}|_{x_{1}=z_{1},
x_{2}=z_{2}}
\end{equation}
defined in the region $|z_{2}|>|z_{1}|>0$ 
 are analytic extensions of each other.
In the definition of l.g.r. conformal 
intertwining algebra, skew-symmetry 
can be replaced by  commutativity.\epf
\end{prop}

In \cite{H6}, generalized rationality for intertwining operator 
algebras is also proved.   In this paper, we do not need the 
full generalized rationality. What we need is a generalization below
of Lemma 4.1 in \cite{H6} to
l.g.r. conformal 
intertwining algebras. Its proof is 
completely the same as the proof of that lemma in \cite{H6}:

\begin{prop}\label{correl1}
Let $(W, 
\mathcal{A}, \{\mathcal{V}_{a_{1}a_{2}}^{a_{3}}\}, {\bf 1}, \omega)$ 
be a l.g.r. conformal 
intertwining algebra. For any $a_{1}, a_{2}, a_{3}, a_{4},
a_{5}\in \mathcal{A}$, any $\mathcal{Y}_{1}\in \mathcal{
V}_{a_{1}a_{5}}^{a_{4}}$ and $\mathcal{Y}_{2}\in\mathcal{
V}_{a_{2}a_{3}}^{a_{5}}$, there exist a 
multivalued analytic function defined on $M^{2}=\{(z_{1},
z_{2})\in \mathbb{C}^{2}\;|\;z_{1}, z_{2}\ne 0, z_{1}\ne
z_{2}\}$,
$\mathcal{Y}^{a, i}_{3}
\in \mathcal{
V}_{a_{1}a_{2}}^{a}$,  $\mathcal{Y}^{a,i}_{4}\in \mathcal{
V}_{aa_{3}}^{a_{4}}$, $\mathcal{Y}^{a,i}_{5}
\in \mathcal{
V}_{a_{2}a}^{a_{4}}$ and $\mathcal{Y}^{a, i}_{6}\in \mathcal{
V}_{a_{2}a_{3}}^{a}$, for $a\in \mathcal{A}$ and $i=1, \dots, k$,
 such that (\ref{prod}), (\ref{iter}) and (\ref{op-prod}) 
are restrictions of this function to their domains 
$|z_{1}|>|z_{2}|>0$, $|z_{2}|>|z_{1}-z_{2}|>0$ and
$|z_{2}|>|z_{1}|>0$, respectively. \epf
\end{prop}

\renewcommand{\theequation}{\thesection.\arabic{equation}}
\renewcommand{\thethm}{\thesection.\arabic{thm}}
\setcounter{equation}{0}
\setcounter{thm}{0}

\section{Complete and completely-extendable conformal intertwining 
algebras}

In this section, we first introduce complete and completely-extendable
conformal intertwining (super)algebras. In particular, we have 
complete and completely-extendable
intertwining operator (super)algebras.
Then we give a construction
of complete conformal intertwining (super)algebras generalizing 
the construction of intertwining operator (super)algebras given in 
\cite{H1}, 
\cite{H4} and 
\cite{H6}. (The details of this generalization will be given in 
another paper \cite{H8}.) We then give examples.

\begin{defn}
{\rm A {\it complete conformal intertwining (super)algebra}
is a l.g.r. conformal intertwining (super)algebra 
$$(W, 
\mathcal{A}, \{\mathcal{V}_{a_{1}a_{2}}^{a_{3}}\}, {\bf 1}, 
\omega)$$
satisfying the following conditions:

\begin{enumerate}

\item Every irreducible 
$W^{e}$-module is completely reducible.

\item The set $\{W^{a}\}_{a\in \mathcal{A}}$
is a  complete set 
of representatives of
equivalence classes of
irreducible $W^{e}$-modules.

\item For any $a_{1}, a_{2}, a_{3}\in \mathcal{A}$, 
$\mathcal{V}_{a_{1}a_{2}}^{a_{3}}$ is the spaces of 
all intertwining operators of type ${W^{a_{3}}\choose
W^{a_{1}}W^{a_{2}}}$.

\end{enumerate}

A {\it completely-extendable conformal intertwining (super)algebra} 
is  a l.g.r. conformal intertwining (super)algebra
$$(W, 
\mathcal{A}, \{\mathcal{V}_{a_{1}a_{2}}^{a_{3}}\}, {\bf 1}, 
\omega)$$
satisfying the following conditions:

\begin{enumerate}

\item Every irreducible 
$W^{e}$-module is completely reducible.

\item The $W^{e}$-modules $W^{a}$, $a\in \mathcal{A}$,
are irreducible.

\item For a complete set 
$\{W^{a^{\#}}\}_{a^{\#}\in \mathcal{A}^{\#}}$
of representatives  of
equivalence classes of
irreducible $W^{e}$-modules, the direct sum $W^{\#}=\coprod_{a^{\#}\in 
\mathcal{A}^{\#}}W^{a^{\#}}$
equipped with the index set $\mathcal{A}^{\#}$,
the vacuum and 
the Virasoro element of $W^{e}$ and the spaces 
$(\mathcal{V}^{\#})_{a_{1}^{\#}a_{2}^{\#}}^{a_{3}^{\#}}$,
$a_{1}^{\#}, a_{2}^{\#}, a_{3}^{\#}\in \mathcal{A}^{\#}$, of all 
intertwining operators of type ${a_{3}^{\#}\choose a_{1}^{\#}a_{2}^{\#}}$, 
is also a l.g.r. conformal
intertwining (super)algebra (which 
by definition is complete). 

\end{enumerate}

By definition, different 
l.g.r. conformal
intertwining (super)algebras obtained from different sets of representatives
of equivalence classes of
irreducible $W^{e}$-modules are naturally isomorphic. We 
 call $W^{\#}$ a
{\it complete extension} of $W$.

{\it Complete conformal vertex
(super)algebras}, {\it completely-extendable conformal vertex
(super)algebras}, {\it complete
vertex operator (super)algebras}, {\it completely-extendable vertex
operator (super)algebras}, {\it meromorphic complete intertwining
algebras} and {\it $\mathbb{Z}$-graded complete intertwining
algebras}
 are defined in the 
obvious way. }
\end{defn}

From the definition, we see that a completely-extendable conformal 
intertwining (super)algebra
is a complete extension of itself if and only if it is complete.

\begin{rema}\label{super2}
{\rm As remarked in Remark \ref{super1}, a complete 
conformal intertwining 
superalgebra can be viewed as a conformal intertwining algebra. 
But it is not complete when viewed as a conformal intertwining algebra.}
\end{rema}

In \cite{H1}, for a rational vertex operator algebra 
satisfying suitable conditions, the associativity of
intertwining operators 
was proved. Since the associativity is the only nontrivial 
property one needs to verify when one tries to construct 
an intertwining operator algebras from representations of 
a rational vertex operator algebra, this result in \cite{H1}
in fact gives a construction of intertwining operator algebras. 
By construction, 
these intertwining operator algebras 
are all complete. The same construction, except
for the change of some signs, works for rational 
vertex operator superalgebras satisfying suitable conditions (see 
the relevant discussions in \cite{HM1}).

In the case of familiar vertex operator (super)algebras (certain 
vertex
operator (super)algebras associated to the Virasoro algebra, affine Lie
algebras, the $N=1$ and $N=2$ Neveu-Schwarz algebras and related algebras), 
the conditions needed were verified in \cite{H3},
\cite{HL4}, \cite{HM1} and \cite{HM2} and
thus we obtain complete intertwining operator (super)algebras from
representations of these vertex operator (super)algebras.

This construction of complete intertwining operator (super)algebras
has a generalization to a construction of complete conformal intertwining 
(super)algebras. Here we state the theorem but its proof will be 
given in a long paper \cite{H8} in preparation on intertwining algebras.

We shall formulate the results only for conformal intertwining algebras.
The super case is the same except for the change of some signs.
To formulate the result precisely, we need 
the following notion which 
in the case of 
vertex operator algebra was first introduced in \cite{H1}:

\begin{defn}
{\rm Let $V$ be a l.g.r. conformal
vertex algebra. We say that products of
intertwining operators
for $V$ satisfy the {\it
convergence and extension property} 
if for any intertwining operators
$\mathcal{Y}_{1}$ and $\mathcal{Y}_{2}$ of types
${W_{0}}\choose {W_{1}W_{4}}$ and ${W_{4}}\choose
{W_{2}W_{3}}$, respectively,
there exists 
an integer $N$
(depending only on $\mathcal{Y}_{1}$ and $\mathcal{Y}_{2}$), and 
for any $w_{(1)}\in W_{1}$,
$w_{(2)}\in W_{2}$, $w_{(3)}\in W_{3}$, $w'_{(0)}\in W'_{0}$, there exist
$j\in \mathbb{N}$, $r_{i}, s_{i}\in \mathbb{R}$, $i=1, \dots, j$, and analytic 
functions $f_{i}(z)$ on $|z|<1$, $i=1, \dots, j$, 
such that
$$
\langle w'_{(0)}, \mathcal{Y}_{1}(w_{(1)}, x_{1})
\mathcal{Y}_{2}(w_{(2)}, x_{2})w_{(3)}\rangle
\lbar_{x_{1}= z_{1}, \;x_{2}=z_{2}}
$$
is absolutely 
convergent when $|z_{1}|>|z_{2}|>0$ and can be analytically extended to  
the multivalued analytic function
$$
\sum_{i=1}^{j}z_{2}^{r_{i}}(z_{1}-z_{2})^{s_{i}}
f_{i}\left(\frac{z_{1}-z_{2}}{z_{2}}\right)
$$
when $|z_{2}|>|z_{1}-z_{2}|>0$. In the case that 
$V$ is a vertex operator algebra, we  require in addition that
\begin{equation}\label{si}
\wt w_{(1)}+\wt w_{(2)}+s_{i}>N,\;\;\;i=1, \dots, j.
\end{equation}}
\end{defn}

\begin{rema}
{\rm It is easy to see that 
 if the associativity of intertwining operators holds
 (see the associativity axiom in Definition \ref{ioa}), then products 
of intertwining operators do satisfy the convergence and extension property.
So this condition is in fact necessary for the  
associativity of intertwining operators. }
\end{rema}

\begin{thm}\label{cons-ioa}
Let $V$ be a l.g.r. conformal
vertex algebra and 
$\{W^{a}\}_{a\in \mathcal{A}}$
a complete set of representatives of  equivalence classes of
irreducible $V$-modules. Assume that every $V$-module is completely
reducible, every irreducible $V$-module
is $\mathbb{R}$-graded, products of
intertwining operators
for $V$ satisfy the
convergence and extension property, and for any 
$a_{1}, a_{2}\in \mathcal{A}$, there are only 
finitely many $a_{3}\in \mathcal{A}$ such that there exist
nonzero intertwining operators of type ${a_{3}\choose a_{1}a_{2}}$. 
Then the 
$\mathbb{R}$-graded vector space
$\coprod_{a\in \mathcal{A}}W^{a}$ 
together with the index set $\mathcal{A}$ with
the special element $e$ such that $W^{e}$ is isomorphic to 
a subalgebra of $V$,
the spaces $\mathcal{V}_{a_{1}a_{2}}^{a_{3}}$, $a_{1}, a_{2}, a_{3}\in
\mathcal{A}$, of all intertwining operators of type
${a_{3}\choose a_{1}a_{2}}$, 
the vacuum ${\bf 1}\in W^{e}$ and the Virasoro element
$\omega\in W^{e}$ is a conformal intertwining algebra. 
If $V$ is in addition
simple in the sense that $V$ as a $V$-module is irreducible,
this conformal intertwining algebra is complete.\epf
\end{thm}

As is mentioned above, the proof of this theorem will be given in 
\cite{H8}.

The construction of intertwining operator algebras obtained in \cite{H1} 
and \cite{H3} can be obtained as an 
easy consequence of the result above. 
Here we state one theorem obtained in \cite{H1} 
and \cite{H3} and formulated in terms of intertwining operator algebras
 first in \cite{H4}. 
We need several generalizations of the notion of module
for a vertex operator algebra. A {\it generalized
$V$-module} is a $\mathbb{C}$-graded vector space equipped with a vertex
operator map satisfying all the axioms for a $V$-module except for the two
grading-restriction conditions. A generalized $V$-module 
is said to 
be {\it locally-grading-restricted} or simply {\it 
l.g.r.} if it
is a $V$-module when $V$ is viewed as a l.g.r. 
conformal vertex algebra. A l.g.r. generalized $V$-module is
said to be {\it lower-truncated} 
if there exists $N\in \mathbb{Z}$ such that
its homogeneous subspace of weight $n$ is
$0$ when the real part of $n$ is less than $N$.

\begin{cor}
Let $V$ be a rational 
vertex operator algebra and 
$\{W^{i}\}_{i=1}^{m}$ 
a set of representatives of equivalence classes of
irreducible $V$-modules. Assume that every irreducible $V$-module
is $\mathbb{R}$-graded, products of
intertwining operators
for $V$ satisfy the
convergence and extension property (including the 
additional condition (\ref{si})), and  
every finitely-generated
lower-truncated  l.g.r. generalized 
$V$-module is a $V$-module. Then the 
$\mathbb{R}$-graded vector space
$\coprod_{i=1}^{m}W^{i}$ together with the finite set 
$\mathcal{A}=\{1, \dots, m\}$ with
the special element $e$ such that $W^{e}$ is isomorphic to a
subalgebra of $V$, 
the spaces $\mathcal{V}_{a_{1}a_{2}}^{a_{3}}$, $a_{1}, a_{2}, a_{3}\in
\mathcal{A}$, of all intertwining operators of type
${a_{3}\choose a_{1}a_{2}}$, 
the vacuum ${\bf 1}\in V$ and the Virasoro element
$\omega\in V$ is a  intertwining operator algebra. If $V$ is in addition 
simple in the sense that $V$ is irreducible as a $V$-module, 
then $\coprod_{i=1}^{m}W^{i}$ is complete.\epf
\end{cor}

The proof of this corollary using Theorem \ref{cons-ioa} will also be 
given in \cite{H8}. 

Now we give examples of complete and complete-extendable conformal 
intertwining operator algebras. Since these examples are mostly highly 
nontrivial, we refer the reader to the references for details except
for algebras given in Examples \ref{moonshine} and \ref{bosons}.

\begin{expl}\label{moonshine}
{\rm Rational vertex operator algebras with themselves as the only
irreducible modules up to isomorphisms: 
Any such algebra is a complete vertex operator
algebra. In particular, the moonshine module $V^{\natural}$ is a
complete vertex operator algebra.}
\end{expl}

\begin{expl}
{\rm Lattice abelian intertwining algebras: In \cite{DL2}, 
Dong and Lepowsky constructed abelian
intertwining algebras associated to nondegenerate lattices. These
abelian intertwining algebras are completely-extendable
conformal intertwining algebras.}
\end{expl}

\begin{expl}\label{bosons}
{\rm Free bosons: Let $\mathfrak{h}_{\mathbb{R}}$ be an $m$-dimensional 
vector space over
$\mathbb{R}$ equipped with a nondegenerate bilinear form $(\cdot,
\cdot)$. We choose a basis $\{h_{i}\}_{i=1}^{m}$
of $\mathfrak{h}_{\mathbb{R}}$ which is orthonormal
in the sense that $(h_{i}, h_{j})=\eta_{ij}$
 where $\eta_{ij}=0$ when $i\ne j$ and $\eta_{ii}=\pm 1$.
Let $\mathfrak{h}=\mathfrak{h}_{\mathbb{R}}
\otimes_{\mathbb{R}}\mathbb{C}$ be the complexification
of $\mathfrak{h}_{\mathbb{R}}$ with the bilinear
form obtained by extending linearly from the 
one on $\mathfrak{h}_{\mathbb{R}}$.  For convenience, we
use the same notations $h_{i}$ to denote
$h_{i}\otimes 1\in \mathfrak{h}$ for $i=1, \dots, 
m$. We view $\mathfrak{h}$ as an abelian Lie algebra and 
consider the $\mathbb{Z}$-graded
untwisted affine Lie algebra $\tilde\mathfrak{h}=\coprod_{n\in \mathbb{Z}}
\mathfrak{h}\otimes t^{n}\oplus \mathbb{C}k\oplus \mathbb{C}d$, its
Heisenberg subalgebra $\hat\mathfrak{h}_{\mathbb{Z}}=\coprod_{n\in \mathbb{Z},
n\ne 0}\mathfrak{h}
\otimes t^{n}\oplus \mathbb{C}k$ and the subalgebra 
$\hat\mathfrak{h}_{\mathbb{Z}}^{-}=\coprod_{n<0}\mathfrak{h}\otimes t^{n}$.  The
symmetric algebra $S(\hat\mathfrak{h}_{\mathbb{Z}}^{-})$ over 
$\hat\mathfrak{h}_{\mathbb{Z}}^{-}$ is a $\mathbb{Z}$-graded 
$\hat\mathfrak{h}_{\mathbb{Z}}$-irreducible $\tilde\mathfrak{h}$-module. 
 There is a unique vertex operator
algebra structure on $S(\hat\mathfrak{h}_{\mathbb{Z}}^{-})$ such that 
$1$ is the vacuum, $\omega=\frac{1}{2}\sum_{i=1}^{m}\eta_{ii}h_{i}(-1)^{2}$
is the Virasoro element and for any $h\in \mathfrak{h}$, $Y(h(-1), x)
=\sum_{n\in \mathbb{Z}}h(n)x^{-n-1}$.

Given any $g\in \mathfrak{h}_{\mathbb{R}}$, 
$S(\hat\mathfrak{h}_{\mathbb{Z}}^{-})\otimes 
\mathbb{C}g$ is a $\hat\mathfrak{h}_{\mathbb{Z}}$-irreducible 
$\tilde\mathfrak{h}$-module such that for any $h\in \mathfrak{h}$,
$h(0)(1\otimes g)=(h, g)(1\otimes g)$. This $\tilde\mathfrak{h}$-module
has a unique irreducible module structure for the vertex operator algebra
$S(\hat\mathfrak{h}_{\mathbb{Z}}^{-})$ such that 
for any $h\in \mathfrak{h}$, $Y(h(-1), x)
=\sum_{n\in \mathbb{Z}}h(n)x^{-n-1}$. It is easy to show that 
any irreducible module for $S(\hat\mathfrak{h}_{\mathbb{Z}}^{-})$
is isomorphic to such a module, such modules with different 
$g\in \mathfrak{h}_{\mathbb{R}}$ are not isomorphic, and any module for 
$S(\hat\mathfrak{h}_{\mathbb{Z}}^{-})$ is completely reducible. 
It is easy to show that $\coprod_{g\in \mathfrak{h}_{\mathbb{R}}}
S(\hat\mathfrak{h}_{\mathbb{Z}}^{-})\otimes 
\mathbb{C}g$ is a l.g.r. conformal intertwining algebra, and by 
definition, it is complete.}
\end{expl}

\begin{expl}
{\rm The $\mathbb{Z}_{2}$-orbifold theory underlying the moonshine
module: In \cite{H2}, the author constructed an abelian intertwining 
algebra on the direct sum of the Leech lattice vertex operator algebra
and its irreducible twisted module. It is a complete intertwining
operator algebra.}
\end{expl}

\begin{expl}
{\rm Minimal models: In \cite{H3}, the author constructed the
intertwining operator algebras
associated to the minimal models. They are complete.}
\end{expl}

\begin{expl}
{\rm WNZW models: In \cite{HL6}, Lepowsky and the author constructed the
intertwining operator algebras associated to the WNZW models. They are 
complete.}
\end{expl}

\begin{expl}
{\rm $N=1$ superconformal minimal models: In \cite{HM1},
Milas and the author constructed the intertwining operator superalgebras
associated to the $N=1$ superconformal minimal models. They are complete.}
\end{expl}

\begin{expl}
{\rm $N=2$ superconformal unitary models: In \cite{HM2},
Milas and the author constructed the intertwining operator superalgebras
associated to the $N=2$ superconformal unitary models. They are complete.}
\end{expl}

\begin{expl}
{\rm Other examples: 
Let $V$ be a rational vertex operator algebra
containing a subalgebra isomorphic to a tensor product algebra
of the vertex operator algebras associated to the minimal models
or the WNZW models or $N=1$ superconformal minimal models or
$N=2$ superconformal unitary models. In \cite{H3} or \cite{HL6} or \cite{HM1} 
or \cite{HM2}, respectively, it was proved that the direct sum of all 
irreducible 
$V$-modules in a complete set of equivalence classes of all 
irreducible modules has a structure of intertwining operator algebra.
It is complete.}
\end{expl}

\renewcommand{\theequation}{\thesection.\arabic{equation}}
\renewcommand{\thethm}{\thesection.\arabic{thm}}
\setcounter{equation}{0}
\setcounter{thm}{0}
\section{Duals of completely-extendable intertwining operator algebras}

In this section, we  construct 
duals of completely-extendable conformal
intertwining algebras and prove their
basic properties. We also introduce nondegenerate 
completely-extendable conformal
intertwining algebras
and prove that for these algebras, their double duals
are equal to themselves.

Let 
$$(W, 
\mathcal{A}, \{\mathcal{V}_{a_{1}a_{2}}^{a_{3}}\}, {\bf 1}, 
\omega)$$
be a completely-extendable conformal intertwining algebra. In this 
and next section, we fix a complete extension
$$(W^{\#}, 
\mathcal{A}^{\#}, \{(\mathcal{V}^{\#})_{a^{\#}_{1}a^{\#}_{2}}
^{a^{\#}_{3}}\}, {\bf 1}, 
\omega)$$
$W$. 
Consider $a^{\#}\in \mathcal{A}^{\#}$
such that for any $a\in \mathcal{A}$, 
$a^{\#}_{1}\in \mathcal{A}^{\#}$ and
any $\mathcal{Y}\in \mathcal{
V}_{aa^{\#}}^{a^{\#}_{1}}$, the image of
$\mathcal{Y}: W^{a}\otimes W^{a^{\#}}\mapsto W^{a^{\#}_{1}}\{x\}$
is in $W^{a^{\#}_{1}}[[x^{-1}, x]]$ (and thus in $W^{a^{\#}_{1}}((x))$).
We denote the set of all such elements $a^{\#}\in \mathcal{A}^{\#}$
by $\mathcal{A}^{\circ}$. We shall use $a^{\circ}, a_{1}^{\circ}, \dots$,
to denote elements of $\mathcal{A}^{\circ}$.
Let 
$W^{\circ}=\coprod_{a^{\circ}\in \mathcal{A}^{\circ}}
W^{a^{\circ}}$. 
First we have the following result:

\begin{prop}\label{intertwining}
Let $a^{\circ}_{1}, a^{\circ}_{2}\in \mathcal{A}^{\circ}$
and  $a^{\#}_{3}\in \mathcal{A}^{\#}$. If
$(\mathcal{V}^{\#})_{a^{\#}_{1}a^{\#}_{2}}^{a^{\#}_{3}}\ne 0$, 
then $a^{\#}_{3}\in
\mathcal{A}^{\circ}$.
\end{prop}
\pf
Suppose that $a^{\#}_{3}\not\in \mathcal{A}^{\circ}$. Then by definition,
there exist $a\in \mathcal{A}$, $a^{\#}_{4}\in \mathcal{A}^{\#}$,
$w_{(a)}\in W^{a}$, $w_{(a^{\#}_{1})}\in W^{a^{\#}_{1}}$,
$w_{(a^{\#}_{2})}\in W^{a^{\#}_{2}}$,
$\mathcal{Y}_{1}\in (\mathcal{V}^{\#})_{aa^{\#}_{3}}^{a_{4}^{\#}}$
and $\mathcal{Y}_{2}\in 
(\mathcal{V}^{\#})_{a^{\#}_{1}a^{\#}_{2}}^{a^{\#}_{3}}$
such that 
$$\mathcal{Y}_{1}(w_{(a)}, x_{1})\mathcal{Y}_{2}(w_{(a^{\#}_{1})}, x_{2})
w_{(a^{\#}_{2})}$$
has terms in nonintegral powers of $x_{1}$.
On the other hand, by Proposition \ref{correl1}, for any 
$w'_{(a^{\#}_{4})}\in (W^{a^{\#}_{4}})'$, 
$$\langle w'_{(a^{\#}_{4})}, 
\mathcal{Y}_{1}(w_{(a)}, x_{1})\mathcal{Y}_{2}(w_{(a^{\#}_{1})}, x_{2})
w_{(a^{\#}_{2})}\rangle_{W^{a^{\#}}}|_{x_{1}=z_{1},
x_{2}=z_{2}}$$
is absolutely convergent in the region $|z_{1}|>|z_{2}|>0$ and 
can be analytically extended to a (multivalued) analytic function $f$
on $M^{2}=\{(z_{1},
z_{2})\in \mathbb{C}^{2}\;|\;z_{1}, z_{2}\ne 0, z_{1}\ne
z_{2}\}$. 
Since 
$$\mathcal{Y}_{1}(w_{(a)}, x_{1})\mathcal{Y}_{2}
(w_{(a^{\#}_{1})}, x_{2})
w_{(a^{\#}_{2})}$$
has nonintegral power terms in $x_{1}$, there exists  
$w'_{(a^{\#}_{4})}\in (W^{a^{\#}_{4}})'$ such that 
$f$ is multivalued in $z_{1}$. But by Proposition \ref{correl1}
again, there exist $\mathcal{Y}^{a^{\#}_{5}, i}_{3}
\in (\mathcal{V}^{\#})_{aa^{\#}_{1}}^{a^{\#}_{5}}$,  
$\mathcal{Y}^{a^{\#}_{5},i}_{4}\in 
(\mathcal{V}^{\#})_{a^{\#}_{5}a^{\#}_{2}}^{a^{\#}_{4}}$, 
$\mathcal{Y}^{a^{\#}_{5},i}_{5}
\in (\mathcal{V}^{\#})_{a^{\#}_{2}a^{\#}_{5}}^{a^{\#}_{4}}$ and 
$\mathcal{Y}^{a^{\#}_{5}, i}_{6}\in 
(\mathcal{V}^{\#})_{aa^{\#}_{2}}^{a^{\#}_{5}}$, 
for $a^{\#}_{5}\in \mathcal{A}^{\#}$ and $i=1, \dots, k$, such that
$f$ is equal to 
$$\sum_{a^{\#}_{5}\in \mathcal{A}^{\#}}\sum_{i=1}^{k}
\langle w_{(a^{\#}_{4})}', \mathcal{Y}^{a^{\#}_{5}, i}_{4}
(\mathcal{Y}^{a^{\#}_{5}, i}_{3}(w_{(a)},
x_{0})w_{(a^{\#}_{1})}, x_{2})w_{(a^{\#}_{2})}\rangle_{W^{a^{\#}_{4}}}
|_{x_{0}=z_{1}-z_{2},
x_{2}=z_{2}}$$
in the region $|z_{2}|>|z_{1}-z_{2}|>0$ and is equal to
and
$$\sum_{a^{\#}_{5}\in \mathcal{A}^{\#}}
\sum_{i=1}^{k}\langle w'_{(a^{\#}_{4})}, 
\mathcal{Y}_{5}^{a^{\#},i}(w_{(a^{\#}_{1})}, x_{2})
\mathcal{Y}_{6}^{a^{\#},i}(w_{(a)}, x_{1})
w_{(a^{\#}_{2})}\rangle_{W^{a^{\#}_{4}}}|_{x_{1}=z_{1},
x_{2}=z_{2}}$$
in the region $|z_{2}|>|z_{1}|>0$. 
Since $a^{\circ}_{1}, a^{\circ}_{2}\in \mathcal{A}^{\circ}$. So
$\mathcal{Y}^{a^{\#}_{5}, i}_{3}(w_{(a)},
x_{0})w_{(a^{\#}_{1})}$ and 
$\mathcal{Y}_{6}^{a^{\#},i}(w_{(a)}, x_{1})
w_{(a^{\#}_{2})}$, $a^{\#}_{5}\in \mathcal{A}^{\#}$,
$i=1, \dots, k$,
have only terms in integral power of $x_{0}$ and $x_{1}$, respectively.
Thus for fixed $z_{2}$, $f$ 
as a function of $z_{1}$ is single-valued near the singularities
$z_{1}=z_{2}$ and $z_{1}=0$. But for fixed $z_{2}$, $f$ 
as a function of $z_{1}$ has only three singularities, $z_{1}=\infty$,
$z_{1}=z_{2}$ and $z_{1}=0$. It is impossible to have such 
an analytic function
which is single-valued near two singularities but multivalued near 
the other singularity. Contradiction. 
\epfv

The main result of this section is the following:

\begin{thm}
The direct sum $W^{\circ}=\coprod_{a^{\circ}\in \mathcal{A}^{\circ}}
W^{a^{\circ}}$, equipped with $\mathcal{A}^{\circ}$, 
the spaces 
$$(\mathcal{V}^{\circ})_{a_{1}^{\circ}a_{2}^{\circ}}^{a_{3}^{\circ}}=
(\mathcal{V}^{\#})_{a_{1}^{\circ}a_{2}^{\circ}}^{a_{3}^{\circ}}$$
of intertwining operators,
the vacuum and the Virasoro element of $W^{e}$, is an 
conformal intertwining algebra. In addition, it is completely extendable 
and its complete extension is equal to the complete extension of 
$W$. 
\end{thm}
\pf
We need only prove the associativity and the skew-symmetry. Let 
$\mathcal{Y}_{1}\in (\mathcal{
V}^{\circ})_{a^{\circ}_{1}a^{\circ}_{5}}^{a^{\circ}_{4}}$ 
and $\mathcal{Y}_{2}\in(\mathcal{
V}^{\circ})_{a^{\circ}_{2}a^{\circ}_{3}}^{a^{\circ}_{5}}$. 
By the associativity for the completely extension $W^{\#}$ of $W$,
there exist $\mathcal{Y}^{a^{\#} i}_{3}
\in (\mathcal{V}^{\#})_{a^{\circ}_{1}a^{\circ}_{2}}^{a^{\#}}$ 
and $\mathcal{Y}^{a^{\#},i}_{4}\in 
(\mathcal{V}^{\#})_{a^{\#}a^{\circ}_{3}}^{a^{\circ}_{4}}$ for all $a^{\#}\in 
\mathcal{A}^{\#}$ and $i=1, \dots, k$
such that the (multivalued) analytic function 
$$
\langle w_{(a^{\circ}_{4})}', 
\mathcal{Y}_{1}(w_{(a^{\circ}_{1})}, x_{1})\mathcal{Y}_{2}(w_{(a^{\circ}_{2})}, 
x_{2})w_{(a^{\circ}_{3})}\rangle_{W^{a_{4}^{\circ}}}|_{x_{1}=z_{1},
x_{2}=z_{2}}
$$
defined in the region
$|z_{1}|>|z_{2}|>0$ 
and the (multivalued) analytic function
\begin{equation}\label{iter2}
\sum_{a^{\#}\in \mathcal{A}^{\#}}\sum_{i=1}^{k}
\langle w_{(a^{\circ}_{4})}', \mathcal{Y}^{a^{\#}, i}_{4}
(\mathcal{Y}^{a^{\#}, i}_{3}(w_{(a^{\circ}_{1})},
x_{0})w_{(a^{\circ}_{2})}, x_{2})w_{(a^{\circ}_{3})}
\rangle_{W^{a^{\circ}_{4}}}
|_{x_{0}=z_{1}-z_{2},
x_{2}=z_{2}}
\end{equation}
defined in the region
$|z_{2}|>|z_{1}-z_{2}|>0$ are equal in the intersection
$|z_{1}|> |z_{2}|>|z_{1}-z_{2}|>0$.
By Proposition \ref{intertwining}, we see that 
$\mathcal{Y}^{a^{\#}, i}_{3}$ 
must be $0$ if $a^{\#}\not \in \mathcal{A}^{\circ}$, proving the 
associativity. 

The skew-symmetry is obvious from the definition.
\epf

The completely-extendable conformal intertwining algebra
$W^{\circ}$ is called the {\it dual} of $W$. If $W^{\circ}=W$,
we say that $W$ is {\it self-dual}.

\begin{prop}\label{***=*}
Let $W$ be a completely-extendable conformal intertwining  algebra.
Then $\mathcal{A}\subset 
(\mathcal{A}^{\circ})^{\circ}$,
$W\subset (W^{\circ})^{\circ}$, $((\mathcal{A}^{\circ})^{\circ})^{\circ}=
\mathcal{A}$ and
$((W^{\circ})^{\circ})^{\circ}=W^{\circ}$.
\end{prop}
\pf
Let $a\in \mathcal{A}$.
Then by the definition of $W^{\circ}$, for any 
$a^{\circ}\in \mathcal{A}^{\circ}$, $a^{\#}\mathcal{A}$ and
any $\mathcal{Y}\in \mathcal{
V}_{aa^{\circ}}^{a^{\#}}$, the image of
$\mathcal{Y}: W^{a}\otimes W^{a^{\circ}}\mapsto W^{a^{\#}}\{x\}$
is in $W^{a^{\#}_{1}}[[x^{-1}, x]]$. Equivalently,
$\Omega_{-1}(\mathcal{Y}): W^{a^{\circ}}\otimes W^{a}\mapsto W^{a^{\#}}\{x\}$
is in $W^{a^{\#}_{1}}[[x^{-1}, x]]$. Since $\Omega_{-1}$ is a linear 
isomorphism, we conclude that for any $a^{\circ}\in \mathcal{A}^{\circ}$,
$a^{\#}\in \mathcal{A}$ and
any $\mathcal{Y}\in \mathcal{
V}_{a^{\circ}a}^{a^{\#}}$, the image of 
$\mathcal{Y}: W^{a^{\circ}}\otimes W^{a}\mapsto W^{a^{\#}}\{x\}$
is in $W^{a^{\#}}[[x^{-1}, x]]$. Thus by definition, $a\in 
(\mathcal{A}^{\circ})^{\circ}$. Since $a$ is arbitrary, $\mathcal{A}\subset 
(\mathcal{A}^{\circ})^{\circ}$ and $W\subset
(W^{\circ})^{\circ}$.

Now we know that $\mathcal{A}\subset 
((\mathcal{A}^{\circ})^{\circ})^{\circ}$ and $W^{\circ}\subset
((W^{\circ})^{\circ})^{\circ}$. To prove the second conclusion, 
we need only prove that $((\mathcal{A}^{\circ})^{\circ})^{\circ}\subset 
\mathcal{A}^{\circ}$. Let $a^{\circ\circ\circ}\in 
((\mathcal{A}^{\circ})^{\circ})^{\circ}$. By definition, 
for any $a^{\circ\circ}\in (\mathcal{A}^{\circ})^{\circ}$, 
$a^{\#}\in \mathcal{A}$ and
any $\mathcal{Y}\in \mathcal{
V}_{a^{\circ\circ}a^{\circ\circ\circ}}^{a^{\#}}$, the image of
$\mathcal{Y}: W^{a^{\circ\circ}}\otimes W^{a^{\circ\circ\circ}}
\mapsto W^{a^{\#}}\{x\}$
is in $W^{a^{\#}}[[x^{-1}, x]]$. Since $\mathcal{A}\subset 
(\mathcal{A}^{\circ})^{\circ}$,
we have, in particular, for any $a\in \mathcal{A}$,
$a^{\#}\in \mathcal{A}$ and
any $\mathcal{Y}\in \mathcal{
V}_{aa^{\circ\circ\circ}}^{a^{\#}}$, the image of
$\mathcal{Y}: W^{a}\otimes W^{a^{\circ\circ\circ}}
\mapsto W^{a^{\#}}\{x\}$
is in $W^{a^{\#}}[[x^{-1}, x]]$. Thus by definition, 
$a^{\circ\circ\circ}\in \mathcal{A}^{\circ}$, proving 
$((\mathcal{A}^{\circ})^{\circ})^{\circ}\subset 
\mathcal{A}^{\circ}$. 
\epfv

\begin{expl}
{\rm Any vertex operator algebra in Example \ref{moonshine} 
is {\it self-dual}.}
\end{expl}

\begin{expl}
{\rm Any simple l.g.r. conformal vertex algebra satisfying the conditions 
in Theorem \ref{cons-ioa} is completely extendable by the theorem.
Its dual is its complete extension.}
\end{expl}

We need the following notion:

\begin{defn}
{\rm A completely-extendable conformal intertwining algebra
$$(W, 
\mathcal{A}, \{\mathcal{V}_{a_{1}a_{2}}^{a_{3}}\}, {\bf 1}, 
\omega)$$
 is said to be {\it nondegenerate}
if it is the dual of a completely-extendable conformal intertwining 
algebra. It is said to be {\it degenerate} if it is not
nondegenerate.}
\end{defn}

For simplicity, we shall often call a nondegenerate 
completely-extendable conformal intertwining 
algebra simply a {\it nondegenerate conformal intertwining 
algebra}.

We have:

\begin{thm}
Let $W$ be a completely-extendable conformal intertwining algebra.
Then $W$ is nondegenerate if and only if $(W^{\circ})^{\circ}=W$.
\end{thm}
\pf
The ``if'' part follows from the definition. The only if part follows
from Proposition \ref{***=*}.
\epfv

The examples of complete-extendable
conformal intertwining  algebras given
in the preceding section are all nondegenerate.
Here we give an example of degenerate  completely-extendable 
intertwining operator algebra.

\begin{expl}
{\rm Let $W^{0}$ be a vertex operator algebra satisfying the conditions in 
Theorem \ref{cons-ioa} and $W^{1}$ an irreducible $W^{0}$-module. 
We consider $W=W^{0}\oplus W^{1}$ and $\mathcal{A}=\{0, 1\}$. 
Let $\mathcal{V}_{00}^{0}$  be the one-dimensional space spanned by the 
vertex operator defining the vertex operator algebra $W^{0}$,
$\mathcal{V}_{01}^{1}$  the one-dimensional space spanned by the 
vertex operator defining the $W^{0}$-module $W^{1}$, 
$\mathcal{V}_{10}^{1}$ the one-dimensional space spanned by the 
intertwining operator obtained from the 
vertex operator defining the $W^{0}$-module $W^{1}$ using skew-symmetry,
and let $\mathcal{V}_{00}^{1}$, $\mathcal{V}_{01}^{0}$, 
$\mathcal{V}_{10}^{0}$,
$\mathcal{V}_{11}^{0}$ and $\mathcal{V}_{11}^{1}$ be $0$.
Let $\mathbf{1}$ and $\omega$ be the vacuum and Virasoro element of 
$W^{0}$.
Then it is easy to verify that $(W, \mathcal{A}, 
\{\mathcal{V}_{a_{1}a_{2}}^{a_{3}}\}, \mathbf{1}, \omega)$ 
is an intertwining operator algebra. Since 
$W^{0}$  satisfies the conditions in 
Theorem \ref{cons-ioa}, $W$ is completely extendable. 
Now choose $W^{0}$ and $W^{1}$ such that $(W^{0})'$ and $(W^{1})'$ are
 isomorphic to $W^{0}$ and
$W^{1}$, respectively, as $W^{0}$-modules. Then 
the space of intertwining operators of type
${(W^{0})'\choose W^{1}(W^{1})'}$ is isomorphic to the space of 
intertwining operators of type
${W^{0}\choose W^{1}W^{1}}$.
But we know that there is a nonzero intertwining operator of type
${(W^{0})'\choose W^{1}(W^{1})'}$ contragredient to 
the intertwining operator spanning the space $\mathcal{V}^{1}_{10}$
(see \cite{HL2}). Thus the space of all intertwining operators of type
${(W^{0})'\choose W^{1}(W^{1})'}$ is nonzero, and consequently
the space of 
intertwining operators of type
${W^{0}\choose W^{1}W^{1}}$ is nonzero. So we see that for the intertwining
operator algebra $W$, $\mathcal{V}_{11}^{0}=0$ is not equal to the
space of all intertwining operators of type ${W^{0}\choose W^{1}W^{1}}$. 
By definition, $W$ is degenerate.}
\end{expl}

\renewcommand{\theequation}{\thesection.\arabic{equation}}
\renewcommand{\thethm}{\thesection.\arabic{thm}}
\setcounter{equation}{0}
\setcounter{thm}{0}

\section{Codes, lattices and 
completely-extendable intertwining operator algebras}

In this section, for the convenience of the reader, we first recall
definitions and properties of linear binary codes and lattices.  Then we
give the correspondence among codes, lattices and
nondegenerate conformal intertwining  algebras using a table. 
A large part of the
correspondence was known to 
or was conjectured by Frenkel, Lepowsky, Meurman and
Goddard. The main  new thing in this section is that the correspondence 
given here
is mathematically precise because of the 
precise conformal-field-theoretic analogues of codes and lattices
and the results obtained in the preceding sections.

A ({\it binary linear}) {\it code of length $n\in \mathbb{N}$} is a subspace
of the vector space $\mathbb{Z}_{2}^{n}$ over $\mathbb{Z}_{2}$. The
dimension of the subspace is called the {\it dimension} of the
code. Elements of a code are called {\it codewords}. Any element $S$
of $\mathbb{Z}_{2}^{n}$ is a linear combination of the basis elements
$(1, 0, \dots, 0), \dots, (0, \dots, 0, 1)$.  The number of nonzero
coefficients in the linear combination is called the {\it weight} of
$S$ and is denoted by $|S|$. A code $\mathcal{S}$ is said to be {\it
even} if $|S|\in 2\mathbb{Z}$ for all $S\in \mathcal{S}$ and {\it
doubly even} if $|S|\in 4\mathbb{Z}$ for all $S\in \mathcal{S}$. The
usual dot product on $\mathbb{Z}_{2}^{n}$ gives a natural nondegenerate
symmetric bilinear form on $\mathbb{Z}_{2}^{n}$.  The orthogonal space of a
code $\mathcal{S}$ in $\mathbb{Z}_{2}^{n}$ with respect to this
bilinear form is again a code. It is called the {\it dual code} of
$\mathcal{S}$ and is denoted $\mathcal{S}^{\circ}$. A code is called
{\it self-dual} if it is equal to its dual code.
For any code $\mathcal{S}$, we have a polynomial
$$W_{\mathcal{S}}(q)=\sum_{S\in \mathcal{S}}q^{\swt S}.$$

{\small \begin{table}
\begin{center}
\begin{tabular}{|p{.8in}|p{1.7in}|p{2.3in}|}\hline
Codes &Lattices &Completely-extendable conformal 
intertwining algebras \\\hline
Even codes&Integral lattices&Meromorphic completely-extendable
conformal intertwining algebras\\ \hline
Doubly even codes&Even lattices&$\mathbb{Z}$-graded completely-extendable
conformal intertwining algebras\\ \hline
&Positive-definite lattices&Completely-extendable intertwining operator 
algebras\\ \hline
&Positive definite and integral lattices&Meromorphic completely-extendable
intertwining operator algebras\\ \hline
&Positive definite and even lattices&$\mathbb{Z}$-graded 
completely-extendable 
intertwining operator algebras\\ \hline
Golay code&Leech lattice&Moonshine module
vertex operator algebra\\ \hline
Lengths&Ranks&Central charges (or ranks)\\\hline
Weights&Square lengths&Weights\\\hline
Dual codes&Dual lattices&Duals of completely-extendable 
conformal 
intertwining algebras\\\hline
Self-dual codes&Self-dual lattices &Self-dual completely-extendable 
conformal 
intertwining algebras\\\hline
&Nondegenerate lattices &Nondegenerate completely-extendable 
conformal 
intertwining algebras\\\hline
$\mathbb{F}_{2}^{n}$&$L_{\mathbb{Q}}$&Complete extensions\\\hline
$W_{\mathcal{S}}(q)$&$\Theta_{L}(q)$&$\chi_{W}(q)$\\\hline
\end{tabular}
\end{center}
\end{table}}

A {\it (rational) lattice of rank} $n\in \mathbb{N}$ is a rank $n$
free abelian group $L$ equipped with a rational-valued symmetric
$\mathbb{Z}$-bilinear form $(\cdot, \cdot)$. A lattice is {\it
nondegenerate} if its form is nondegenerate. 
Let $L$ be a lattice.  For $m\in
\mathbb{Q}$, we set $L_{m}=\{\alpha\in L\;|\;( \alpha, \alpha
)=m\}$. The lattice $L$ is said to be {\it even} if $L_{m}=0$ for any
$m\in \mathbb{Q}$ which is not an even integer. A lattice $L$ is said to be
{\it integral} if the form is integral valued and to be {\it positive
definite} if the form is positive definite. Even lattices are
integral. Let $L_{\mathbb{Q}} =L\otimes_{\mathbb{Z}} \mathbb{Q}$. Then
$L_{\mathbb{Q}}$ is an $n$-dimensional vector space over $\mathbb{Q}$
in which $L$ is embedded and the form on $L$ is extended to a
symmetric $\mathbb{Q}$-bilinear form on $L_{\mathbb{Q}}$, still denoted
by $(\cdot, \cdot)$. The lattice
is nondegenerate if and only if this form on $L_{\mathbb{Q}}$ is
nondegenerate.  The {\it dual} of $L$ is the set
$L^{\circ}=\{\alpha\in L_{\mathbb{Q}}\; |\;( \alpha, L)\subset
\mathbb{Z}\}$. This set is a lattice if and only if $L$ is nondegenerate,
and in this case, $L^{\circ}$ has as a basis the dual basis of a given
basis.  The lattice $L$ is said to be {\it
self-dual} if $L=L^{\circ}$. This is equivalent to $L$ being integral
and {\it unimodular}, which means that $|\det((\alpha_{i},
\alpha_{j}))|=1$.   For any lattice $L$, we have the
associated theta series
$$\Theta_{L}(q)=\sum_{\alpha\in L}q^{\frac{1}{2}(\alpha, \alpha)}.$$

Given a completely-extendable 
conformal 
intertwining algebra $W$ such that $\dim W_{(n)}<\infty$ 
for $n\in \mathbb{C}$, let 
$$\chi_{W}(q)=q^{-\frac{c}{24}}
\sum_{n\in \mathbb{C}}(\dim W_{(n)})q^{n}.$$
It is called the {\it character} or {\it graded dimension} or 
the {\it partition function} of $W$.

Now see the table for the precise correspondence among codes, lattices
and completely-extendable conformal intertwining algebras. This table
fills some blanks in a table 
in \cite{G} and gives more items than the one in
\cite{G}. There are constructions of lattices from codes and abelian
intertwining algebras from lattices (see \cite{CN} and
\cite{DL2}). 

In fact this correspondence can be formulated more conceptually in
terms of the language of categories and related results can be
established. Some of the constructions mentioned above should be
viewed as canonical functors and others should be viewed as
compositions of these canonical functors with some twist functors in
the categories of lattices and completely-extendable conformal
intertwining algebras. We shall discuss this formulation and these
results in a future paper.

{\small \sc Department of Mathematics, Rutgers University,
110 Frelinghuysen Rd., Piscataway, NJ 08854-8019}

{\em E-mail address}: yzhuang@math.rutgers.edu


\begin{thebibliography}{FLM2}


\bibitem[B]{B}
R.~E.~Borcherds,
Vertex algebras, Kac-Moody algebras, and the Monster,
{\em Proc. Natl. Acad. Sci. USA} {\bf 83} (1986), 3068--3071.

\bibitem[C]{C}
J. H. Conway,
A characterisation of Leech's lattice, 
{\em Invent. Math.} {\bf 7} (1969), 137--142. 

\bibitem[CN]{CN}
J. H. Conway, N. J. Sloane, 
{\em Sphere packings, lattices and groups},
3rd edition, Grundlehren der Mathematischen Wissenschaften, Vol. 290, 
Springer-Verlag, New York, 1999. 


\bibitem[DGM]{DGM}
L. Dolan, P. Goddard and P. Montague, Conformal field theory of
twisted vertex operators, {\em Nucl. Phys.} {\bf B338} (1990),
529--601.

\bibitem[DL1]{DL1}
C.~Dong and J.~Lepowsky, Abelian intertwining algebras --- a generalization
of vertex operator algebras, in:
{\em Proc. AMS Summer Research Institute on Algebraic Groups
and Their Generalizations, Pennsylvania State University, 1991,}  ed.
W. J. Haboush and B. J. Parshall, Proc. Symp. Pure Math., 
Vol 56, Part 2, American Mathematical Society, Providence,
1994, 261--293.

\bibitem[DL2]{DL2}
C.~Dong and J.~Lepowsky, {\em Generalized Vertex Algebras and Relative
Vertex Operators}, Progress in Math., Vol. 112, Birkh\"{a}user,
Boston, 1993.

\bibitem[FFR]{FFR}
A. J. Feingold, I. B. Frenkel, J. F. X. Ries, {\em
Spinor construction of vertex
operator algebras, triality and $E_{8}^{(1)}$}, Contemporary Math.,
Vol. 121, Amer. Math. Soc., Providence, 1991.

\bibitem[FHL]{FHL}
I.~B. Frenkel, Y.-Z. Huang and J.~Lepowsky,
On axiomatic approaches to vertex operator algebras and modules,
preprint, 1989;
{\em Memoirs Amer. Math. Soc.} {\bf 104}, 1993.


\bibitem[FLM1]{FLM1}
I.~B. Frenkel, J.~Lepowsky and A.~Meurman,
A natural representation of the Fischer-Griess Monster with
the modular function $J$ as character,
{\em Proc. Natl. Acad. Sci. USA} {\bf 81} (1984),
3256--3260.

\bibitem[FLM2]{FLM}
I.~B. Frenkel, J.~Lepowsky, and A.~Meurman,
{\em Vertex operator algebras and the Monster},
Pure and Appl. Math., Vol. 134, Academic Press, New York, 1988.

\bibitem[G]{G}
P. Goddard, Meromorphic conformal field theory,
in: {\em Infinite-dimensional Lie algebras and
groups (Luminy-Marseille, 1988)}, ed. V. Kac, Adv. Ser. Math. Phys., 
Vol. 7, World Sci. Publishing, Teaneck,  1989, 556--587.


\bibitem[H1]{H0}
Y.-Z. Huang,
Geometric interpretation of vertex operator algebras,
{\em Proc. Natl. Acad. Sci. USA} {\bf 88} (1991), 9964--9968.


\bibitem[H2]{H1}
Y.-Z. Huang, A theory of tensor products for module categories
for a vertex operator algebra, IV, {\em J. Pure Appl. Alg.} {\bf 100}
(1995), 173-216.

\bibitem[H3]{H2}
Y.-Z. Huang, A nonmeromorphic extension of the moonshine module 
vertex operator algebra, in: {\em Moonshine, the Monster
and related topics, Proc. Joint Summer Research Conference, 
Mount Holyoke, 1994,} ed. C. Dong and G. Mason, 
Contemporary Math., Vol. 193, Amer. Math. Soc., Providence, 1996, 
123--148.

\bibitem[H4]{H3}
Y.-Z. Huang, Virasoro vertex operator algebras,
(nonmeromorphic) operator product expansion and the tensor product
theory, {\em J. Alg.} {\bf 182} (1996), 201--234.

\bibitem[H5]{H4}
Y.-Z. Huang, Intertwining operator algebras, genus-zero modular
functors and genus-zero conformal field theories,  in: {\em Operads:
Proceedings of Renaissance Conferences}, ed. J.-L. Loday,
J. Stasheff, and A. A. Voronov, Contemporary Math., Vol. 202,
Amer. Math. Soc., Providence, 1997,  335--355.

\bibitem[H6]{H4.5}
Y.-Z. Huang, 
{\em Two-dimensional conformal geometry and vertex operator algebras}, 
Progress in Mathematics, Vol. 148, 1997,
Birkh\"{a}user, Boston.


\bibitem[H7]{H5}
Y.-Z. Huang, Genus-zero modular functors and intertwining operator algebras,
{\em Internat. J. Math.} {\bf 9} (1998), 845--863.


\bibitem[H8]{H6}
Y.-Z. Huang, Generalized rationality and a
``Jacobi identity'' for intertwining operator algebras,
{\em Selecta
Mathematica, New Series}, to appear.

\bibitem[H9]{H8}
Y.-Z. Huang, Intertwining algebras, in preparation. 



\bibitem[HL1]{HL1}
Y.-Z. Huang and J. Lepowsky, A theory of tensor products for module
categories for a vertex operator algebra, I, {\em Selecta
Mathematica, New Series} {\bf 1} (1995), 699-756.

\bibitem[HL2]{HL2}
Y.-Z. Huang and J. Lepowsky, A theory of tensor products for module
categories for a vertex operator algebra, II, {\em Selecta
Mathematica, New Series} {\bf 1} (1995), 757--786.

\bibitem[HL3]{HL3}
Y.-Z. Huang and J. Lepowsky, Tensor products of modules for a vertex
operator algebras and vertex tensor categories, in:
     {\em Lie Theory and Geometry,
in honor of Bertram Kostant,}
ed. R. Brylinski, J.-L. Brylinski, V. Guillemin, V. Kac,
Birkh\"{a}user, Boston, 1994, 349--383.

\bibitem[HL4]{HL4}
Y.-Z. Huang and J. Lepowsky, A theory of tensor products for module
categories for a vertex operator algebra, III, {\em J. Pure
Appl. Alg.} {\bf 100} (1995),  141-171.

\bibitem[HL6]{HL6}
Y.-Z. Huang and J. Lepowsky, Intertwining operator algebras and vertex
tensor categories for
affine Lie algebras,  {\it Duke Math. J.} {\bf 99} (1999), 113--134.


\bibitem[HM1]{HM1}
Y.-Z. Huang and A. Milas, Intertwining operator superalgebras and 
vertex tensor categories for superconformal algebras, I, to appear.

\bibitem[HM2]{HM2}
Y.-Z. Huang and A. Milas, Intertwining operator superalgebras and 
vertex tensor categories for superconformal algebras, II, {\em
Trans. Amer. Math. Soc.}, to appear.

\end{thebibliography}
\end{document}